\documentclass{article}
\usepackage{amssymb}
\usepackage{amsmath}
\usepackage{graphicx}

\newtheorem{theorem}{Theorem}

\newtheorem{definition}[theorem]{Definition}

\newtheorem{lemma}[theorem]{Lemma}

\begin{document}
\title{Matching and geometry of 2-step nilpotent Lie groups}

\author {Babak Hasanzadeh\\
Department of Mathematics, Faculty of Science\\
Azarbaijan Shahid Madani University\\
Tabriz, Iran\\
E-mail: babakmath@gmail.com}

\maketitle


\renewcommand{\thefootnote}{}


\renewcommand{\thefootnote}{\arabic{footnote}}
\setcounter{footnote}{0}
\begin{abstract}
 In this paper we study contact structure on 2-step nilpotent, Heisenberg type  Lie groups. We decompose this Lie groups to center and orthogonal complement, then investigate properties of both orthogonal Lie subgroups. Finally, we provide a connection between matchings in groups and field extensions and 2-step nilpotent Lie groups.

\end{abstract}

\section{Introduction and preliminaries}
2-step nilpotent Lie groups play an important role in mathematics and have a special significance. The geometry of nonabelian nilpotent Lie group studied by many authors but our main resources are P. Eberlein' s papers on 2-step nilpotent Lie groups in \cite{tk} and \cite{tkp}.  
In this section we bring up basic definitions and theorems, and we will do some calculations for the following sections. Throughout the paper Lie groups are assumed to be simply connected. A Lie group H of a Lie group G is a subgroup which is also a submanifold and $ \mathfrak{g} $ is Lie algebra of G, for any $ X,Y \in \mathfrak{g} $, we have
\begin{align*}
[X,Y]=\nabla_{X}Y-\nabla_{Y}X.
\end{align*}
We define the linear map $ adX:\mathfrak{g} \longrightarrow \mathfrak{g}$ by
 $Y \longmapsto [X,Y] \nonumber $.

Let G be a Lie group equipped by Riemannian left invariant metric, by using only these identities and combining a few permutations of variables obtain the formula. 
\begin{align}
g(\nabla_{X}Y,Z) &=\dfrac{1}{2}\lbrace Xg(Y,Z)+Yg(X,Z)- Zg(X,Y) \\   
                         &+g([X,Y],Z)+g([Z,X],Y)-g([Y,Z],X)\rbrace, \nonumber
\end{align}
 for each nonzero left invariant vector fields   $ X,Y,Z \in \mathfrak{g} $. $ \nabla $ is covariant derivative in this paper first three terms are vanish.
 If (G,g) is a Lie group equipped by left invariant metric, ad is skew adjoint if
\begin{align*}
g(adX (Y),Z)=-g(Y,adX (Z)).
\end{align*}
 If $ X \in Z( \mathfrak{g}) $ and $ Y,Z \in  \mathfrak{g} $ ,we have
\begin{align*}
X g(Y,Z)=g(\nabla_{Y+Z}X,Y).
\end{align*}
In the next definition, we define the most basic concept that will be use in this paper.
\begin{definition}
A nilpotent Lie group is a Lie group G which is connected and whose Lie algebra is nilpotent Lie algebra $ \mathfrak{g} $, that is, it's Lie algebra has a sequence of ideals of $ \mathfrak{g} $ by $ \mathfrak{g}^{0}=\mathfrak{g} $, $ \mathfrak{g}^{1}=[\mathfrak{g},\mathfrak{g}] $,$ \mathfrak{g}^{2}=[\mathfrak{g},\mathfrak{g}^{1}] $,..., $ \mathfrak{g}^{i}=[\mathfrak{g},\mathfrak{g}^{i-1} ]  $. Also $ \mathfrak{g} $ is called nilpotent if $ \mathfrak{g}^{n}=0 $ for some n.\cite{tg}\\ 
\end{definition}
The following subgroup
\begin{equation*}
Z(G)=\lbrace x \in G:xy=yx, \forall y \in G \rbrace,
\end{equation*}
is called the center of G, it is a Lie subgroup with corresponding Lie subalgebra
\begin{equation*}
Z(\mathfrak{g}) =\lbrace X \in \mathfrak{g} : [X,Y]=0, \forall Y\in \mathfrak{g} \rbrace.
\end{equation*}
Let $ Z^{\perp}(\mathfrak{g}) $ is orthogonal complement of $ Z(\mathfrak{g}) $ with left invariant metric g.
Each element $ Z \in Z(\mathfrak{g}) $ defines a skew symmetric linear map $ j(Z) :Z^{\perp}(\mathfrak{g}) \rightarrow Z^{\perp}(\mathfrak{g})  $ given by $ j(Z)X=(adX)^{\ast}$  for all  $ X \in Z^{\perp}(\mathfrak{g}) $, where $ (adX)^{\ast} $ is the adjoint of adX relative to g, equally and more usefully $ j(Z) $ is defined by the equation 
\begin{align}
g(j(Z)X,Y)=g([X,Y],Z), 
\end{align}
for all $ X,Y \in Z^{\perp}(\mathfrak{g}) $.\\
Following definition is the most practical concept in present paper.
\begin{definition}
 A finite dimensional Lie algebra $ \mathfrak{g} $ is 2-step nilpotent if $ \mathfrak{g} $ is not abelian and $ [\mathfrak{g},[\mathfrak{g},\mathfrak{g}]]=0 $. A Lie group G is 2-step nilpotent if its Lie algebra $ \mathfrak{g} $ is 2-step nilpotent. In the other words, a Lie algebra  $ \mathfrak{g} $ is 2-step nilpotent if $ [\mathfrak{g},\mathfrak{g}] $ is non-zero and lay in the center of $\mathfrak{g}$.\cite{tk}\\
\end{definition}
\begin{theorem}
Let G be a 2-step nilpotent Lie group of Heisenberg type. Then 
\begin{align*}
[X,j(Z)X]=\mid X \mid^{2} Z
\end{align*}
for all elements $ X \in Z^{\perp}(\mathfrak{g}) $ and $ Z \in Z(\mathfrak{g}) $.\cite{tk}
\end{theorem}
\begin{lemma}
Let G be a 2-step nilpotent Lie group of Heisenberg type. Then 
\begin{align*}  
 [X_{n},j(Z_{c})X_{n}]=\mid X_{n} \mid^{2}Z_{c},
\end{align*} 
 for all elements $ X,Z \in \mathfrak{g} $. \cite{tk}\\
\end{lemma}
Because G is a Heisenberg Lie group from \cite{tk} we know
 \begin{align*} 
g(j(Z)X,j(Z^{'})X)=g(Z,Z^{'})\mid X \mid^{2}.
\end{align*}
Some formulas and definitions are introduced, that will be used throughout this paper.
We may identify an element of $ \mathfrak{g} $ with a left invariant vector field on G since $ T_{e}G $ may be identified with $ \mathfrak{g} $. If X,Y are left invariant vector fields on G, then $ \nabla_{X}Y $ is left invariant as well. One has the following formulas:
\begin{flushleft}
(a)$\nabla_{X}Y=\dfrac{1}{2}[X,Y],$  \\
(b)$\nabla_{X}Z=\nabla_{Z}X =-\dfrac{1}{2}j(Z)X,$  \\
(c)$\nabla_{Z}Z_{1}=0,$\\ 
\end{flushleft}
 for all $X,Y \in Z^{\perp}(\mathfrak{g}) $ and $ Z,Z_{1} \in Z(\mathfrak{g}) $\cite{tk}.\\ 
Next definitions are important concepts in present peper.
\begin{definition}
 A 2-step nilpotent Lie algebra $ \mathfrak{g} $ is nonsingular if $ adX:\mathfrak{g} \rightarrow Z(\mathfrak{g}) $ is surjective for each $ X \in Z^{\perp}(\mathfrak{g}) $. A  2-step nilpotent Lie group G  is nonsingular if its Lie algebra $ \mathfrak{g} $ is nonsingular.\cite{tk}\\ 
\end{definition}
\begin{definition}
A totally geodesic subgroup of G is a connected Lie subgroup N, such that N is a totally geodesic Lie subgroup as a submanifold of G. A subalgebra $ \mathfrak{n} $ of $ \mathfrak{g} $ is totally geodesic if $ \nabla_{X}Y \in  \mathfrak{n} $ whenever $ X,Y \in  \mathfrak{n} $.\cite{tk}\\
\end{definition}
Let G denote a 2-step nilpotent Lie group with a left invariant metric g and $ \mathfrak{g} $ denote the Lie algebra of G. Write $ \mathfrak{g}=Z(\mathfrak{g}) \oplus  Z^{\perp}(\mathfrak{g}) $ where $  Z^{\perp}(\mathfrak{g}) $ its orthogonal complement of center $ Z(\mathfrak{g}) $, hence for any $ X \in \mathfrak{g} $ we have $ X=X_{c}+X_{n} $, $ X_{c} \in Z(\mathfrak{g}) $ and $ X_{n} \in  Z^{\perp}(\mathfrak{g}) $. For any $ X,Y \in \mathfrak{g} $, by straightforward calculations we obtain following results and complete this section.
\begin{align*}
&(i)\nabla_{X_{n}}Y_{n}=\dfrac{1}{2}[X,Y], \\
&(ii)\nabla_{X}X=-j(X_{c})X_{n}, \\
&(iii)\nabla_{X}Y=\dfrac{1}{2}[X_{n},Y_{n}]-\dfrac{1}{2}j(Y_{c})X_{n}-\dfrac{1}{2}j(X_{c})Y_{n}.\\
\end{align*}
The following definition is given from \cite{cm}, see also \cite{co},\cite{cn} for more results on matching. 
\begin{definition}
A 2-step nilpotent Lie group is said to have the matching property if for any $ A,B \subseteq G $, $ \mid A\mid =\mid B\mid $, there exist a matching from A to B. 
\end{definition}

\section{contact structure and Heisenberg Lie groups} 
In this section we first introduce almost contact structure, and defined this structure on 2-step nilpotent Lie group, then we try to find new properties of these  Lie groups as a smooth manifolds.\\  
 Let $ M $ be an odd dimensional smooth manifold with a left invariant metric $ g $. Denote by $ TM $ the Lie algebra of vector fields on $ M $. Then $ M $ is said to be an almost contact metric manifold if there exists a tensor $ \phi $ of type $ (1,1) $, a vector field $ \xi $ called structure vector field and $ \eta $, the dual 1-form of $ \xi $ satisfying the following
\begin{align}
 \phi^{2}X=-X+\eta(X)\xi ,    g(X,\xi)=\eta(X) 
\end{align}
\begin{align}
\eta(\xi)=1    ,    \phi(\xi)=0     ,   \eta\circ\phi=0 
\end{align}
\begin{align}
 g(\phi X,\phi Y)=g(X,Y)-\eta(X)\eta(Y), 
\end{align}
for any $ X,Y \in TM $. In this case  
\begin{align*}
 g(\phi X,Y)= - g(X,\phi Y), 
\end{align*}
the fundamental 2-form $ \Phi $ on $ M $ is given by
\begin{align*}
\Phi (X,Y)=g(X,\phi Y)
\end{align*}
and the manifold is said to be contact metric manifold if $ \Phi=d\eta $. Every Heisenberg Lie group is contact metric Lie group.\\
{\bf Example.} Let $ \lbrace X_{1},X_{2},X_{3},Z_{1},Z_{2} \rbrace $ be basis of a 5-dimensional real vector space $\mathfrak{g}$, equip with a Lie algebra structure defined by the bracket relations
\begin{align*}
[X_{1},X_{2}]=Z_{1}, \quad [X_{2},X_{3}]=Z_{1}, \quad [X_{1},X_{3}]=Z_{2}.
\end{align*} 
It's trivial $\mathfrak{g}$ is 2-step nilpotent and nonsingular, with left invariant metric g. By easy computation we have
\begin{align*}
j(Z_{1})X_{1}=X_{2}, \quad j(Z_{1})X_{2}=X_{3}, \quad j(Z_{2})X_{1}=X_{3}.
\end{align*} 
Now we can write 
\begin{align*}
g(j(Z_{1})^{2}X_{1},X_{3}) &=-g(j(Z_{1})X_{1},j(Z_{1})X_{3})\\
                                       &=g(X_{2},X_{2}).
\end{align*} 
Therefore $\mathfrak{g}$ is nonsingular, 2-step nilpotent Lie algebra.\\
It is defined Nijenhuis torsion of $ \phi $
\begin{align*}
[\phi ,\phi](X,Y)=\phi^{2}[X,Y]+[\phi X,\phi Y]-\phi [\phi X,Y]-\phi [X,\phi Y],
\end{align*}
and 
\begin{align*}
N^{(1)}(X,Y) &=[\phi ,\phi](X,Y)+2d\eta (X,Y)\xi, \\
N^{(2)}(X,Y) &=(\mathcal{L}_{\phi X}\eta)(Y)-(\mathcal{L}_{\phi Y}\eta)(X),\\
N^{(3)}        &=(\mathcal{L}_{\xi}\phi)X,\\
N^{(4)}        &=(\mathcal{L}_{\xi}\eta)X.\\
\end{align*}
An almost contact structure $ (\phi ,\xi , \eta) $ is normal if and only if these four tensors vanish. \cite{cam}\\
In continue, we study a Lie group as an almost contact manifold and applying the Lie algebra and new results are obtained about the components of almost contact structure and their features. In almost contact structure $ (\phi ,\xi , \eta) $, from $ N^{(3)} $ we have the following theorem.\\
Let N be a subgroup of G and $ \mathfrak{n} $ denote the Lie algebra of N, if $ \mathfrak{n} $ is an Abelian, totally geodesic subalgebra of N, then either $ \mathfrak{n} \subseteq Z(\mathfrak{g}) $ or $ \mathfrak{n} $ is an Abelian subspace of  $ Z^{\perp}(\mathfrak{g})  $, \cite{tkp}. The angle $ \theta (X) $, $ 0 \leq \theta (X) \leq \dfrac{\pi}{2} $, between $ \phi X $ and $ \mathfrak{n} $ is called the Wirtinger angle of X. If the Wirtinger angle $ \theta $ is constant, it's called slant angle. If $  \mathfrak{n} $ be a slant Lie subalgebra, then  H is called slant Lie subgroup, In the other words slant submanifold. If $ \theta =0 $, then N is called invariant subgroup. \cite{mo}\\
Next theorem show important feature of 2-step nilpotent Lie groups.
\begin{theorem}
Let (G,g) be an almost contact, 2-step nilpotent Lie group with a left invariant metric and $ \mathfrak{g} $  denote the Lie algebra of G. Then
\begin{align*}
&a) \phi (Z(\mathfrak{g}))\neq 0, \\
&b) \phi (Z^{\perp}(\mathfrak{g})) \neq 0. \\
\end{align*}
\end{theorem}

Next theorem show that what is the relationship between Abelian subalgebra and invariant subalgebra with j.
\begin{theorem}
 Let G be a 2-step nilpotent, Heisenberg, Lie group with a left invariant metric, and N is a totally geodesic Lie subgroup of G. Let $ \mathfrak{n}  $ and $ \mathfrak{g} $ denote their Lie algebras, respectively. If for any $ X \in \mathfrak{n} $ we have $ j(X_{c})X=X $, then $  \mathfrak{n} $ is an Abelian Lie subalgebra or $ \mathfrak{n} \subset Z^{\perp}(\mathfrak{g}) $.
\end{theorem}
The following theorem gives the necessary sufficient condition for 2-step nilpotent Lie groups having the matching property. See \cite{cp},\cite{ck} for more details. 
\begin{theorem}
A Heisenberg Lie group G passess the matching property if and only if G has no finite proper subgroup. 
\end{theorem}

\end{document}